
\input amstex
\documentstyle{amsppt}
\magnification 1200
\NoRunningHeads
\NoBlackBoxes
\document

\def\Hom{\text{Hom}}

\def\a{\frak a}

\def\Ua{U_q(\tilde\g)}
\def\U2{{\Ua}_2}
\def\g{\frak g}
\def\h{\frak h}

\def\<{\langle}
\def\>{\rangle}
\def\o{\otimes}

\topmatter
\title Quantization of Poisson algebraic groups and Poisson homogeneous spaces
\endtitle
\author {\rm {\bf Pavel Etingof and David Kazhdan} \linebreak
\vskip .1in
Department of Mathematics\linebreak
Harvard University\linebreak
Cambridge, MA 02138, USA\linebreak
e-mail: etingof\@math.harvard.edu\linebreak kazhdan\@math.harvard.edu}
\endauthor
\endtopmatter

\centerline{October 20, 1995}
\vskip .1in

\centerline{\bf Abstract}
\vskip .1in
This paper consists of two parts.
In the first part we show that
any Poisson algebraic group over a field of characteristic zero
and any Poisson Lie group admits a local quantization.
This answers positively a question of Drinfeld
and generalizes the results of \cite{BFGP} and \cite{BP}.
In the second part we apply our techniques of quantization to obtain
some nontrivial examples of quantization of Poisson homogeneous spaces.

\vskip .1in
\centerline{\bf 1. Quantization of Poisson algebraic and Lie groups.}
\vskip .1in

Below we will freely use the notation from our previous paper
\cite{EK}.
\vskip .1in

1.1. {\it Poisson algebras and manifolds.}

Let $k$ be a field of characteristic $0$, $B$ be a commutative algebra
over $k$. A Poisson bracket on $B$ is a Lie bracket
$\{,\}:B\o B\to B$
which satisfies the Leibnitz rule $\{f,gh\}=\{f,g\}h+g\{f,h\}$,
$f,g,h\in B$. A commutative algebra $B$ equipped with a Poisson bracket
is called a Poisson algebra.

Let $B,C$ be Poisson algebras. Then $B\o C$
has a natural structure of a
Poisson algebra, defined by $\{b_1\o c_1,b_2\o c_2\}=b_1b_2\o\{c_1,c_2\}+
\{b_1,b_2\}\o c_1c_2$.

Let $X$ be an algebraic variety over $k$. Denote by $\Cal O_X$
the structure sheaf of $X$. $X$ is called a Poisson algebraic variety
if the sheaf $\Cal O_X$ is equipped with the structure of a sheaf
of Poisson algebras. Similarly one defines the notions of a smooth or
complex analytic Poisson manifold.

If $X,Y$ are Poisson manifolds, then the product $X\times Y$ has a natural
structure of a Poisson manifold. This is a consequence of the fact that
the tensor product of Poisson algebras is a Poisson algebra.
\vskip .1in

1.2. {\it Poisson groups.}
Let $\a$ be a finite-dimensional Lie algebra over $k$.
Then the ring $k[[\a]]$ of formal power series on $\a$ is a commutative
topological Hopf algebra. We will denote by $A$ the ``spectrum'' of $k[[\a]]$
We call $A$ the formal group associated to $\a$, and
say that $k[[\a]]$ is the structure ring of $A$. Of course,
the precise meaning of spectrum has to be clarified,
but this is irrelevant for us,
as we only use the structure ring
$k[[\a]]$.

Let $A$ denote either
a formal or an affine algebraic group over $k$.
When we refer to both cases we will use the term ``Poisson group''.

Let $O_A$ be the structure ring of $A$.
Let $\Delta:O_A\to O_{A\times A}$
be the standard coproduct. In the case of algebraic groups, this coproduct
is given by the formula $\Delta(f)(x,y)=f(xy)$, $x,y\in A$.

Let $\{,\}$ be a
Poisson bracket on $O_A$.
This bracket defines a Poisson bracket on $O_{A\times A}$,
as explained in Section 1.1.
We say that the Poisson bracket $\{,\}$
is compatible with the group structure if the coproduct $\Delta$ is a
homomorphism of Poisson algebras. If $A$ is equipped with a Poisson
bracket compatible with the group structure, it is called a Poisson
formal (respectively, algebraic) group.

Any Poisson algebraic group
canonically defines a Poisson formal group.
Indeed, let $I$ be the ideal of all elements of $O_A$ vanishing at the
identity. It is clear that the Poisson bracket maps
$I^m\times I^n$ to $I^{m+n-2}$, so it extends to the projective limit
of $O_A/I^n$, which is exactly $k[[\a]]$, where $\a$ is the Lie
algebra of $A$. So we obtain the structure of a Poisson
formal group on $k[[\a]]$.
\vskip .1in

1.3. {\it Poisson formal groups and Lie bialgebras.}

\proclaim{Theorem 1.1} (Drinfeld) The category of
Poisson formal groups over $k$ is equivalent to the category
of Lie bialgebras over $k$.
\endproclaim

\demo{Proof} Let $k[[\a]]$ be a Poisson
formal group. This ring carries a formal series topology and
is a commutative topological Hopf algebra. The space $k[[\a]]^*$ of
continuous linear functionals on $k[[\a]]$ is the universal enveloping algebra
$U(\a)$ of $\a$. Therefore, the adjoint map
to $\{,\}$ defines a skew-symmetric cobracket $\delta: U(\a)\to U(\a)\o U(\a)$
satisfying the Jacobi identity, the Leibnitz rule
$$
(\tilde\Delta_0\o 1)(\delta(x))=(1\o\delta)(\tilde\Delta_0(x))+
s_{23}(\delta\o 1)(\tilde\Delta_0(x))
$$
 (where $s_{23}$ is the permutation of
the second and the third component), and the compatibility relation
$\delta(xy)=\delta(x)\tilde\Delta_0(y)+\tilde\Delta_0(x)\delta(y)$.
It is easy to show that $\delta(\a)\subset \a\o\a$, and
that $\delta$ defines a Lie bialgebra structure on $\a$.
Conversely, if $\a$ is a Lie bialgebra, then the cobracket
$\delta$ on $\a$ can be uniquely extended to a map
$\delta: U(\a)\to U(\a)\o U(\a)$ using the compatibility relation.
This map is skew-symmetric and satisfies the Leibnitz rule
and the Jacobi identity, so the adjoint of this map is a
Poisson bracket on $k[[\a]]$ which defines the structure of a
Poisson formal group on $k[[\a]]$. It is easy to check that this correspondence
defines an equivalence of categories.
\enddemo

4. {\it Quantization of Poisson groups.}

\proclaim{Definition 1.1}
Let $A$ be a Poisson formal or algebraic group.
A quantization of $A$ is a product $*:O_A[[h]]\o O_A[[h]]\to O_A[[h]]$
and coproduct
$\Delta: O_A[[h]]\to O_{A\times A}[[h]]$ on
the topological $k[[h]]$-module $O_A[[h]]$,
such that $O_A[[h]],*,\Delta$ is a topological Hopf algebra,
and
$$
f*g=fg+\frac{h}{2}\{f,g\}+O(h^2), \Delta(f)=\Delta_0(f)+O(h^2), f,g\in O_A,
$$
where $\Delta_0$ is the undeformed coproduct in $O_A$.
\endproclaim

\proclaim{Proposition 1.2} Any Poisson formal group admits a quantization.
\endproclaim

\demo{Proof}
Let $\a$ be a Lie bialgebra.
Let $U_h(\a)$ be the quantization of $\a$ introduced in \cite{EK}.
We identify it, as a $k[[h]]$-module, with $U(\a)[[h]]$,
by the map $\mu: U(\a)[[h]]\to U_h(\a)$
defined in Section 4.1 of \cite{EK}. Consider
the space $ H=U_h(\a)^*$ of all continuous $k[[h]]$-linear functionals
$U_h(\a)\to k[[h]]$. This space carries a natural weak topology and
is isomorphic to $k[[\a]][[h]]$ as a topological $k[[h]]$-module.
Moreover, this space carries a natural structure of a topological
Hopf algebra, dual to the structure of a topological Hopf algebra
in $U_h(\a)$. It is clear from the definition of $U_h(\a)$
that $ H/h H$ is isomorphic to $k[[\a]]$ as a topological Hopf algebra,
and for any $f,g\in k[[\a]]\subset  H$ we have the expansions
$f*g=fg+\frac{1}{2}h\{f,g\}+O(h^2), \Delta(f)=\Delta_0(f)+O(h^2)$,
where $*,\Delta$ are the product and coproduct in $ H$, and $\Delta_0$
is the coproduct in $k[[\a]]$. Thus, the algebra $ H$ is a quantization
of the Poisson formal group $k[[\a]]$.
\enddemo
\vskip .1in

1.5. {\it Local quantization.}

\proclaim{Definition 1.2}
Let $H=(O_A[[h]],*,\Delta)$ be a quantization of a Poisson group $A$.
Let $D_A$ be the algebra of differential operators on $A$.
We say that $H$ is a {\it local quantization} of $A$ if
each coefficient of the $h$-expansion
of $f*g$ is of the form
$\sum_{i=1}^mD_ifD_i^\prime g$, where $D_i,D_i^\prime\in D_A$, and
each coefficient of the $h$-expansion
of $\Delta(f)$ is of the form $D\Delta_0(f)$, where $D\in D_{A\times A}$.
\endproclaim

{\bf Remark.}
The notion of a local quantization is important for the following reason.
If we have a local quantization of a Poisson algebraic group $A$,
we can extend this quantization to the structure sheaf $\Cal O_A$ of $A$.
Namely,
for an open set $U\subset A$ denote by $\Cal O_A(U)$ the algebra
of regular functions on $U$. Then
we have an associative product $*:\Cal O_A(U)[[h]]\o \Cal O_A(U)[[h]]
\to \Cal O_A(U)[[h]]$
and a coproduct $\Delta: \Cal O_A(U)[[h]]\to \Cal O_{A\times A}
(m^{-1}(U))[[h]]$,
where $m:A\times A\to A$ is the product in $A$.
These operations are obtained by extension of the formulas for $O_A$
to $\Cal O_A(U)$, which is possible due to the locality of the quantization.
In particular, we obtain a
sheaf of associative algebras on $A$ which is a quantization of the structure
sheaf of $A$ as a sheaf of Poisson algebras.
\vskip .1in

1.6. {\it Existense of local quantization.}

In this section we will prove the following theorem,
which is the main result of Chapter 1:

\proclaim{Theorem 1.3} Any Poisson formal group
$G$ over $k$ admits a local quantization. Any Poisson algebraic group
$G$ over $k$ admits a local quantization.
\endproclaim

To prove Theorem 1.3, we need the following result.

\proclaim{Proposition 1.4} (i)  The quantization $H$ of
$k[[\a]]$ constructed in Section 1.4
 is local.

(ii) If $k[[\a]]$ is the Poisson formal group associated
to a Poisson algebraic group $A$, then
the coefficients of the quantization constructed in section 1.4
are algebraic
differential operators.
More precisely, each coefficient of the $h$-expansion
of $f*g$ is of the form
$\sum_{i=1}^mD_ifD_i^\prime g$, where $D_i,D_i^\prime\in D_A$, and
each coefficient of the $h$-expansion
of $\Delta(f)$ is of the form $D\Delta_0(f)$, where $D\in D_{A\times A}$.
\endproclaim

\demo{Proof}
We first show the locality of the product $*$.
To be consistent with the notation of \cite{EK}, we denote
the Lie bialgebra $\a$ by $\g_+$,
the dual Lie bialgebra to $\g_+$ by $\g_-$, and
the double $\g_+\oplus\g_-$ by $\g$.

In \cite{EK}, for every Drinfeld associator
and any Lie bialgebra $\g_+$ we defined the topological Hopf algebra
$U_h(\g_+)$ which is a quantization of $\g_+$.
The algebra $U_h(\g_+)$ was defined as the space
$\Hom(M_+\o M_-,M_-)$, where $M_+,M_-$ are the Verma modules over $\g$,
introduced in Section 2.3 of \cite{EK}, and $\Hom$ is taken
in the Drinfeld category $\Cal M$ associated to $\g$. Analogously to \cite{EK},
we identify $U_h(\g_+)$ with $M_-[[h]]$ by the map
$\nu_+: U_h(\g_+)\to M_-[[h]]$ given by $\nu_+(f)=f(1_+\o 1_-)$.
Then the coproduct in $U_h(\g_+)$ is given by the following formula
(formula (4.9) in \cite{EK}):
$$
\tilde\Delta(x)=
J^{-1}i_-(x), x\in M_-,
$$
$J\in U(\g)\o U(\g)[[h]]$ is given by formula (3.1) of \cite{EK}.
Because $M_-=U(\g_+)1_-$, and the map $i_-: M_-\to M_-\o M_-$ corresponds
to the undeformed coproduct on $U(\g_+)$,
 we obtain the following formula for the *-product on
$k[[\g_+]]$:
$$
f*g=p((J^{-1})^*(f\o g)),
$$
where $p: k[[\g_+]]\o k[[\g_+]]\to k[[\g_+]]$ denotes the undeformed product.
Thus, the locality property of the star-product is a consequence of the
following lemma.

\proclaim{ Lemma 1.5}
The algebra $U(\g)$ acts in $k[[\g_+]]=M_-^*$ by differential
operators. In addition, if $k[[\g_+]]$ is the formal group corresponding to
a Poisson algebraic group $G_+$, then these differential operators are in
$D_{G_+}$.
\endproclaim

{\it Proof of the Lemma.} It is enough to show that the Lie algebra
$\g$ acts in $k[[\g_+]]$ by vector fields from $D_{G_+}$.
This is obvious for the subalgebra $\g_+\subset \g$, which acts
by left-invariant vector fields. So, we have to prove the statement
for the subalgebra $\g_-$. To do this,
we recall that $J=1+hr/2+O(h^2)$, where $r\in\g_+\o\g_-$ is the
classical $r$-matrix of $\g$. This shows that for any
$f,g\in k[[\g_+]]$ we have $\{g,f\}=p(r(f\o g))=\sum_j x_j^+f\cdot x_j^-g$,
where $x_j^+$ is a basis of $\g_+$, $x_j^-$ is the dual basis of $\g_-$.
Thus, for any $g\in O_{G_+}$ the vector field $\sum_j(x_j^-g)x_j^+$
is from $D_{G_+}$. Since $\{x_j^+\}$ are a basis of the space of left-invariant
vector fields on $G_+$, we see that the functions $x_j^-g$ belong to $O_{G_+}$
for all $j$, i.e. $x_j^-\in D_{G_+}$. The lemma is proved.
\vskip .03in

Now let us show the locality of the coproduct. Let $x,y\in M_-$. Let
$X,Y\in \Hom(M_+\o M_-,M_-)=U_h(\g_+)$ be such that $X(1_+\o 1_-)=x,
Y(1_+\o 1_-)=y$. Let us compute the product $Z=YX$ in $U_h(\g_+)$
and the vector $z=Z1_-\in M_-$. As follows from Chapter 4 of \cite{EK},
the element $z$ is given by the formula
$$
z=X(1\o Y)\Phi^{-1}(1_+\o 1_+\o 1_-),
$$
where $\Phi$ is the associator. Therefore,
it follows from Lemma 2.1 in \cite{EK} that every coefficient of the
$h$-expansion of $z$ is representable as a finite sum of terms of the form
$$
X(1\o Y)(a\o \tilde\Delta_0(b))(1_+\o 1_+\o 1_-), a,b\in U(\g),
$$
where $\tilde\Delta_0$ is the standard coproduct of $U(\g)$.
Using the fact that $Y$ is an intertwiner, we can rewrite this expression as
$$
X(a1_+\o by).
$$
Since $X$ is an intertwiner, we can further rewrite this
expression as a linear combination of expressions of the form
$$
cX(1_+\o dy), c,d\in U(\g).
$$
Let $t\in U(\g_+)$ be such that $t1_-=dy$. Then
$cX(1_+\o dy)=cX\tilde\Delta_0(t)(1_+\o 1_-)=ctx$.

Now consider the elements $\hat x,\hat y\in U(\g_+)$ such that $\hat x1_-=x$,
$\hat y1_-=y$. For any $u\in U(\g)$, $s\in U(\g_+)$ denote by $u(s)$ the
element of $U(\g_+)$ such that $u(s)1_-=us1_-$. We have shown
that each coefficient of the $h$-expansion of $z$ is a linear combination of
terms of the form $c(d(\hat y)\hat x)$, $c,d\in U(\g)$.
Therefore, using Lemma 1, we can conclude that the coefficients
of $\Delta(f)$ for $f\in k[[\g_+]]$ are given by differential
operators from $D_{G_+\times G_+}$ acting on
$\Delta_0(f)$. The proposition is proved.
\enddemo

\demo{Proof of Theorem 1.3} Let $A$ be a Poisson algebraic group whose Lie
algebra is $\a$.
Proposition 1.3 implies that the formulas which define the quantization
of $k[[\a]]$ described in Section 1.4
automatically define a local quantization of $A$.
\enddemo
\vskip .1in

\proclaim{Theorem 1.6} Let $G$ be a Poisson formal or Poisson algebraic group
whose Lie algebra is a quasitriangular Lie bialgebra. Then $G$ admits
a local formal quantization such that the comultiplication is the same as
in the classical case.
\endproclaim

\demo{Proof} The theorem is an immediate consequence of Theorem 1.3 and
Theorem 6.1 of \cite{EK}, which says that a quasitriangular Lie bialgebra
$\g$ admits a quantization $U_h(\g)$ which is isomorphic to
$U(\g)[[h]]$ as a topological algebra.
\enddemo

{\bf Remark.} Deformations with preserved coproduct are called
preferred deformations\cite{BFGP}.
Theorem 1.6 for Poisson Lie groups
was proved in the reductive case
in \cite{BFGP}, then in general by Bidegain and Pinczon
\cite{BP} and P.Cartier (private communication), using
the existence result of \cite{EK}.
\vskip .1in

1.7. {\it Quantization of Poisson Lie groups.}

One may consider the case when $k=\bold R$ or  $\bold C$, and $A$ is a
Lie group over $k$. In this case $O_A$ is the algebra of smooth
(respectively, holomorphic) functions on $A$. The notions of a Poisson
structure on $A$ and its quantization is
introduced analogously to the algebraic case, and we have a result
analogous to Theorem 1.3.

\proclaim{Theorem 1.7} Any Poisson Lie group admits a local quantization.
\endproclaim

The construction of this quantization is the same as in the algebraic case.

8. {\it Functoriality of quantization.}

The universality of the quantization $U_h(\g_+)$, proved in
Chapter 10 of \cite{EK}, implies that
the quantization of Poisson formal, algebraic and Lie
groups obtained above is functorial. In other words,
we have the following result.

\proclaim{Theorem 1.8} Let $\phi: G_1\to G_2$ be a morphism of Poisson
groups. Then the induced map
$\phi^*: O_{G_2}[[h]]\to O_{G_1}[[h]]$ commutes with
the deformed product $*$ and deformed coproduct $\Delta$,
i.e. defines a homomorphism of topological Hopf algebras.
\endproclaim

\vskip .1in
\centerline{\bf 2. Quantization of Poisson homogeneous spaces.}
\vskip .1in

2.1. {\it Poisson $G$-manifolds.}

Let $G$ be a formal or algebraic group over $k$.
Let $\g$ be the Lie algebra of $G$.
By a $G$-manifold we mean a formal (respectively, algebraic) variety
$X$ with a left action of $G$. More specifically,
in the formal case a $G$-manifold $X$
is a commutative topological algebra $O_X=k[[X]]$ isomorphic to the algebra of
formal power series in $\text{ dim}X$ variables, together with a coaction
$\Delta_0: k[[X]]\to k[[\g]]\o k[[X]]$ of
the Hopf algebra $k[[\g]]$ such that $\Delta_0$ is a homomorphism of algebras.
In the algebraic case, a $G$-manifold is a variety $X$
with a sheaf of commutative algebras $O^X$ (the structure sheaf of $X$)
together with a morphism of sheaves of algebras
$\Delta_0: m^*O^X\to O^{G\times X}$, where $m$ is the
action of $G$ on $X$, $m:G\times X\to X$.

Let $G$ be a Poisson group.
Assume that a $G$-manifold $X$
is equipped with a Poisson bracket $\{,\}$.
We say that $X$ equipped with $\{,\}$ is a Poisson $G$-manifold if
$\Delta_0$ is a homomorphism of (sheaves of) Poisson algebras.
In particular, if $G$ acts transitively on $X$, we say that $X$ is
a Poisson homogeneous space.

{\bf Example.} Let $H\subset G$ be a subgroup, $\h$ be the Lie algebra of $H$.
Assume that $\h$ is a coideal in $\g$, i.e. $\delta(\h)\subset \h\o\g\oplus
\g\o\h$. In this (and only this)
case the Poisson bracket on $G$ descends to
the homogeneous space $G/H$. Namely, in the formal case we obtain
a Poisson bracket $\{,\}$ on $k[[\g/\h]]$ which comes from the natural Lie
coalgebra structure on $\g/\h$. In the algebraic case we obtain
a Poisson bracket $\{,\}$ on the algebra $\Cal O_{G/H}(U)$ for any open set
$U\subset G/H$. In both cases the map $\Delta_0$ is a homomorphism
of Poisson algebras, so $G/H$ is a Poisson $G$-manifold.
We call such a manifold a Poisson homogeneous space of group type.

{\bf Remarks.} 1. Not every Poisson homogeneous space
is of group type. For example, if a group $G$ with trivial Poisson structure
acts transitively on a Poisson manifold $X$ by Poisson automorphisms,
then $X$ is a Poisson homogeneous space, but the Poisson structure on
$X$ is not obtained from the Poisson structure on $G$.

2. In the case of algebraic homogeneous spaces
 we use the language of sheaves, since the
homogeneous space $G/H$ may not be affine even if $G,H$ are affine, and so
the ring of globally defined functions on $G/H$ could contain only constants.
\vskip .1in

2.2. {\it Equivariant quantization of Poisson $G$-manifolds.}

Let $G$ be a Poisson group, and $O_{h,G}$ be a local quantization of
the algebra $O_{G}$.
In the algebraic case, let $\Cal O_{h,G}$ be a quantization of the sheaf
of regular functions on $G$.

\proclaim{Definition 2.1} (i)
By an $O_{h,G}$-equivariant quantization of
a Poisson formal $G$-manifold $X$ we mean
the $k[[h]]$-module $O_X[[h]]$ equipped with
a new associative product $*$ and coaction $\Delta: O_X[[h]]\to
O_{h,G}\o O_X[[h]]$, such that $\Delta$ is a homomorphism of
algebras, and
$$
f*g=fg+\frac{h}{2}\{f,g\}+O(h^2), \Delta(f)=\Delta_0(f)+O(h^2),
f,g\in O_X.\ (1)
$$

(ii) By an $\Cal O_{h,G}$-equivariant quantization of
a Poisson algebraic $G$-manifold $X$ we mean
the sheaf of $k[[h]]$-modules $\Cal O_X[[h]]$ equipped with
a new associative product $*$ and homomorphism of sheaves
of algebras on $G\times X$ (with new multiplication), $\Delta:
m^*(\Cal O_X[[h]])\to \Cal O_{G\times X}[[h]]$, such that
equations (1) are satisfied.

(iii)  We say that a quantization of a
Poisson $G$-manifold $X$ is local
if each coefficient of the $h$-expansion
of $f*g$ is of the form
$\sum_{i=1}^mD_ifD_i^\prime g$, $D_i,D_i^\prime\in D_X$, and
each coefficient of the $h$-expansion
of $\Delta(f)$ is of the form $D\Delta_0(f)$, $D\in
D_{G\times X}$, where
$D_X$ is the algebra (sheaf of algebras) of differential
operators on $X$.
\endproclaim

{\bf Remark.} In the formal case,
it is convenient to talk about quantization in terms of dual spaces to
function algebras. Denote
the space
$O_X^*$ of continuous linear functionals on $O_X$
by $T$. Then the structure of a Poisson $G$-manifold on $X$
amounts to an action of $\g$ on $T$ and a structure of a co-Poisson algebra
on $T$ (i.e. a coproduct $\tilde\Delta_0$ and cobracket $\delta$),
such that the map $U(\g)\o T\to T$ defining the action of $\g$ on $T$
is a homomorphism of co-Poisson algebras. Further, the problem
of quantization of $X$
reduces to finding an action of $U_h(\g)$ on $T[[h]]$ and a
$U_h(\g)$-invariant
coassociative coproduct $\tilde\Delta: T[[h]]\to T[[h]]\o T[[h]]$
which has the form
$\tilde\Delta=\tilde\Delta_0+\frac{1}{2}h\delta+O(h^2)$.

{\bf Example.} We say that a Poisson homogeneous space $G/H$ of group type
 is split if
$\h\subset \g$ is a Lie subbialgebra. Let us construct the quantization
of a split homogeneous space. Let $U_h(\g),U_h(\h)$ be the quantizations
of the Lie bialgebras $\g,\h$. The functoriality of quantization implies that
we have a Hopf algebra embedding $U_h(\h)\to U_h(\g)$. Let
$T=Ind_{U_h(\h)}^{U_h(\g)}{\bold 1}$ be a module over $U_h(\g)$.
This module is generated by a vector $1_T$ such that $a1_T=\epsilon(a)1_T$,
$a\in U_h(\h)$, where $\epsilon$ is the counit.
It is clear that $T^*$ is naturally isomorphic to $k[[\g/\h]][[h]]$.
  Let $i_T: T\to T\o T$ be the intertwiner defined by the formula
$i_T(1_T)=1_T\o 1_T$. Define the product on $T^*$
by $f*g=i_T^*(f\o g)$, and the coproduct to be the dual map to the action
$U_h(\g)\o T\to T$. It is easy to check that this defines a
local quantization
of the formal homogeneous space $k[[\g/\h]]$.
Moreover, if the space $G/H$ is algebraic, we automatically obtain
its quantization as well, due to the locality property of the quantization.

{\bf Remark.}
Quantized homogeneous spaces for semisimple Lie groups are considered
in \cite{BFGP}. It is shown there that any preferred deformation of a
compact group descends to a quotient of this group by a normal subgroup.
\vskip .1in

2.3. {\it Quantization of homogeneous spaces arising from Manin quadruples.}

\proclaim{Definition 2.2} A Manin quadruple is a quadruple $(\g,\g_+,\g_-,\h)$,
where $(\g,\g_+,\g_-)$ is a Manin triple, and $\h$ is a Lagrangian Lie
subalgebra in $\g$.
\endproclaim

{\bf Example.} Let $(\g,\g_+,\g_-)$ be a Manin triple, and
$\h_+\subset \g_+$ be a Lie subalgebra which is also a coideal,
i.e. $\delta(\h_+)\subset \h_+\o\g_+\oplus\g_+\o \h_+$. Let
$\h_-$ be the orthogonal complement of $\h_+$ in $\g_-$.
Then $\h_-$ is also a Lie subalgebra and coideal in $\g_-$,
and the direct sum $\h=\h_+\oplus\h_-$ is an isotropic Lie subalgebra in $\g$,
so $(\g,\g_+,\g_-,\h)$ is a Manin quadruple.

Any Manin quadruple defines a formal Poisson
homogeneous space, as follows.
Let $G,H$ denote the formal groups corresponding
to $\g,\h$. For any two functions $f,g$ on $G$ right
invariant under $H$ set
$$
\{f,g\}=\sum_i R_{x_i^+}f R_{x_i^-}g,
$$
where $x_i^+$ is a basis of $\g_+$, $x_i^-$ is the dual basis of $\g_-$, and
$R_xf$ denotes the action of the right invariant vector field corresponding
to $x\in\g$ on a function $f$.

\proclaim{Lemma 2.1} The operation
$\{,\}$ is skew symmetric.
\endproclaim

{\it Proof.} We can identify the space
$k[[\g/\h]]^*$ with the left $U(\g)$-module $T=U(\g)/U(\g)\h$.
Let $1_T$ denote the image of $1$ under the natural projection
$U(\g)\to T$. Then we have
$$
\{f,g\}(a1_T)=(f\o g)(r\tilde\Delta_0(a)(1_T\o 1_T)), a\in U(\g)
$$
where $\tilde\Delta_0$ is the coproduct in $U(\g)$. In particular,
$$
(\{f,g\}+\{g,f\})(a1_T)=(f\o g)((r+r^{op})\tilde\Delta_0(a)(1_T\o 1_T))
$$
The element $\Omega=r+r^{op}$ is $\g$-invariant, and
$\Omega(1_T\o 1_T)=0$ because the algebra $\h$ is Lagrangian.
Thus, $\{f,g\}+\{g,f\}=0$.

\proclaim{Lemma 2.2} The operation $\{,\}$ is a Poisson bracket on $G/H$
which makes $G/H$ into a Poisson homogeneous space
(not necessarily of group type).
\endproclaim

{\it Proof.} The skew-symmetry follows from Lemma 2.1.
The Leibnitz rule is obvious. The Jacobi identity
follows from the classical Yang-Baxter equation for $r$. Consistency
with the Poisson bracket on $G$ is checked directly.

\proclaim{Theorem 2.3}
The Poisson formal homogeneous space $G/H$ admits a local
quantization. Moreover, if $G',H'$ are Poisson algebraic groups
with Lie algebras $\g,\h$ then the quantization of $G/H$ automatically
determines a quantization of $G'/H'$.
\endproclaim

{\it Proof.} Let $T=\text{ Ind}_\h^\g {\bold 1}$ be a $U(\g)$ module.
It is generated by a vector $1_T$ such that $\h 1_T=0$.
We regard $T$ as an object of category $\Cal M$ defined in Section 1.4
of \cite{EK}. Let $i_T: T\to T\o T$ be a morphism in $\Cal M$
defined by the equation $i_T(1_T)=1_T\o 1_T$.

\proclaim{Lemma 2.4} The morphism
$i_T$ is coassociative in $\Cal M$, i.e. $(i_T\o 1)i_T=(1\o i_T)i_T$.
\endproclaim

\demo{Proof}
This statement is proved analogously to Lemma 2.3 in \cite{EK},
using the identity $\Phi(1_T\o 1_T\o 1_T)=1_T\o 1_T\o 1_T$,
which follows from the fact that $\h$ is a Lagrangian subalgebra.
\enddemo

Denote by $F$ the tensor functor from the category $\Cal M$
and the category of representations of $U_h(\g)$, constructed Chapter 2 of
\cite{EK}. Applying this functor to the morphism $i_T$, we obtain a
coassociative coproduct $F(i_T): F(T)\to F(T)\o F(T)$, which is an intertwiner
for $U_h(\g)$.

Now we can define a quantization of the formal homogeneous space $G/H$.
By definition, the space of functions on this quantized homogeneous space
is the $k[[h]]$-module $F(T)^*$. It is clear that $F(T)^*$ is
isomorphic to $O_{G/H}[[h]]$ as a topological $k[[h]]$-module.
 Denote by $*$ the map $F(T)^*\o F(T)^*\to F(T)^*$ dual to the map
$F(i_T)$, and by $\Delta$ the map $F(T)^*\to O_G[[h]]\o F(T)^*$ dual
to the action of $U_h(\g)$ on $F(T)$. Since the product $*$
is associative and $U_h(\g)$-equivariant,
the triple $(F(T)^*,*,\Delta)$ is a quantized homogeneous space.
It is easy to check, analogously to Section 3.4 of \cite{EK},
that $F(T)^*$ is in fact a quantization of the Poisson homogeneous space $G/H$
defined above. This proves Theorem 2.3 in the formal case.
The algebraic case follows immediately
because of the locality of the quantization. The theorem is
proved.

{\bf Remarks.}

1. The results of this Chapter trivially generalize to
smooth and complex analytic homogeneous spaces.
\vskip .1in

2. Not every Poisson homogeneous space can be quantized. For example,
let $X=\bold R^2$ with the standard symplectic structure,
and $G$ be the group of symplectic diffeomorphisms of $X$ with the zero
Poisson-Lie structure. It is clear that $G$ acts transitively on $X$,
and $X$ is a Poisson homogeneous space. Since the Poisson structure on $G$
is zero, the quantization of $X$ as a Poisson homogeneous space
is the same as a $G$-equivariant quantization of $X$ as a symplectic
manifold. It is well known that such a quantization does not exist.

It is an interesting question which Poisson homogeneous spaces can be
quantized. At present we are unable to answer this question.

\centerline{\bf Acknowledgements}
The authors would like to thank A.Alekseev for many inspiring discussions,
P.Cartier for an interesting conversation, and Y.Soibelman
for raising the questions about quantization of homogeneous spaces.

\end